\documentclass[
leqno,          
11pt,           
a4paper         
]{amsart}

\usepackage[colorlinks,citecolor=blue,urlcolor=blue]{hyperref}          
\usepackage{fourier}                                                    
\usepackage{amssymb, amsmath, amsfonts, amsthm}                         
\usepackage{mathtools, cite, enumerate, rotating, framed, lipsum}       
\usepackage{fancybox, bbm, cancel}                                              
\usepackage[dvipsnames]{xcolor}                                         
\usepackage[polish, english]{babel}                                     
\usepackage[T1]{fontenc}                                                
\usepackage[cp1250]{inputenc}                                           

\numberwithin{equation}{section}                        

\newtheoremstyle{mpthm}                         
{0.5\baselineskip}                              
{0.5\baselineskip}                              
{\normalfont\color[named]{MidnightBlue}}        
{0pt}                                           
{\bfseries\color{RubineRed}}                    
{}                                              %
{\newline}
{{\thmname{#1}~}                                %
{\thmnumber{#2}.~}                              %
{\normalfont\color{RubineRed}\thmnote{(#3)}}}   

\newcommand{\thmcount}{equation}                 
\newcounter{specialcounter}

\newtheorem{Thm}[\thmcount]{Theorem}
\newtheorem{Sthm}[specialcounter]{Theorem}
\newtheorem{Cor}[\thmcount]{Corollary}
\newtheorem{Lem}[\thmcount]{Lemma}
\newtheorem{Prop}[\thmcount]{Proposition}
\newtheorem{Rem}[\thmcount]{Remark}
\newtheorem{Defn}[\thmcount]{Definition}
\newtheorem{Ex}[\thmcount]{Example}
\newtheorem{Sex}[specialcounter]{Example}
\newtheorem{Asu}[\thmcount]{Assumption}
\newtheorem{Sol}[\thmcount]{Solution}
\newtheorem*{Thmx}{Theorem}
\newtheorem*{Corx}{Corollary}
\newtheorem*{Lemx}{Lemma}
\newtheorem*{Propx}{Proposition}
\newtheorem*{Remx}{Remark}
\newtheorem*{Defnx}{Definition}
\newtheorem*{Exx}{Example}
\newtheorem*{Asux}{Assumption}
\newtheorem*{Solx}{Solution}

\newcommand \eq[1]{\begin{equation} #1 \end{equation}}
\newcommand \eqx[1]{\begin{equation*}  #1 \end{equation*}}
\newcommand \al[1]{\begin{align} #1 \end{align}}
\newcommand \alx[1]{\begin{align*}  #1 \end{align*}}

\newcommand \spx[1]{\begin{equation*} \begin{split} #1 \end{split} \end{equation*}}
\newcommand \en[1]{\begin{enumerate}  #1 \end{enumerate}}
\newcommand \ite[1]{\begin{itemize}  #1 \end{itemize}}

\newcommand{\sthm}[2]{\begin{Sthm} \label{#1} #2 \end{Sthm}}

\newcommand{\lem}[2]{\begin{Lem} \label{#1} #2 \end{Lem}}

\newcommand{\cor}[2]{\begin{Cor} \label{#1} #2 \end{Cor}}

\newcommand{\prop}[2]{\begin{Prop} \label{#1} #2 \end{Prop}}

\newcommand{\rem}[2]{\begin{Rem} \label{#1} #2 \end{Rem}}

\newcommand{\sex}[2]{\begin{Sex} \label{#1} #2 \end{Sex}}

\newcommand{\pr}[1]{\begin{proof} #1 \end{proof}}
\newcommand{\hidpr}[1]{}

\newcounter{comcount}

\renewcommand{\hline}{\vbox{\hrule width\textwidth height 1pt}\smallskip}

                \newcommand{\e}{\varepsilon}
\newcommand{\w}{\omega}                 
        \newcommand{\la}{\lambda}

  \renewcommand{\aa}{\mathcal{A}}
  
\newcommand{\CC}{\mathbb{C}}

 \newcommand{\Kk}{\mathbf{K}} 
 \newcommand{\Ll}{\mathbf{L}} 
 \newcommand{\Mm}{\mathbf{M}} 
\newcommand{\NN}{\mathbb{N}} \newcommand{\Nn}{\mathbf{N}}

 \newcommand{\Qq}{\mathbf{Q}} \newcommand{\qq}{\mathcal{Q}}
\newcommand{\RR}{\mathbb{R}}  
  
 \newcommand{\Tt}{\mathbf{T}} 
  
  \newcommand{\vv}{\mathcal{V}}

\newcommand{\supp}{\mathrm{supp}}
\newcommand{\8}{\infty}

\renewcommand{\(}{\left(}
\renewcommand{\)}{\right)}

\newcommand{\Rd}{{\RR^d}}

\newcommand{\limj}{\lim_{j\to \8}}

\newcommand{\sumj}{\sum_{j=1}^{\8}}
\newcommand{\sumi}{\sum_{i=1}^{\8}}

\newcommand{\abs}[1]{\left| #1 \right|}
\newcommand{\set}[1]{\left\{ #1 \right\}}
\newcommand{\norm}[1]{\left\| #1 \right\|}
\newcommand{\expr}[1]{\( #1 \)}

\newcommand{\lab}[1]{\label{#1}}

\renewcommand{\LLL}{{L^1(\Rd)}}

\newcommand{\LU}{\L{U}}
\newcommand{\LV}{\L{V}}
\newcommand{\LW}{\L{W}}
\newcommand{\LVW}{\L{V+W}}

\newcommand{\K}[2]{\Kk_{#1}^{#2}}
\renewcommand{\k}[3]{K_{#1}^{#2} (#3)}
\newcommand{\T}[1]{\Tt_{#1}}
\renewcommand{\t}[2]{T_{#1} (#2)}
\renewcommand{\L}[1]{\Ll^{#1}}
\newcommand{\1}{\mathbbm{1}}
\newcommand{\HwQ}{H^1_{at}(\aa_{\w,\qq})}
\newcommand{\HQ}{H^1_{at}(\aa_{\qq})}

\newcommand{\Deltax}{\boldsymbol\Delta}

\title[Atomic decompositions for Hardy spaces related to Schr\"odinger operators]{ Atomic decompositions for Hardy spaces\\ related to Schr\"odinger operators }
\author[ Marcin Preisner ]{ Marcin Preisner }
\address{
Marcin Preisner \newline
\indent Instytut Matematyczny, Uniwersytet Wroc\l awski \newline
\indent pl. Grunwaldzki 2/4, 50-384 Wroc\l aw, Poland }
\email{preisner@math.uni.wroc.pl }

\thanks{ The research was supported by Narodowe Centrum Nauki (NCN) grant nr 2012/05/B/ST1/00672}

\begin{document}
\begin{abstract}
Let $\LU = -\Deltax +U$ be a Schr\"odinger operator on $\Rd$, where $U\in L^1_{loc}(\Rd)$ is a~non-negative potential and $d\geq 3$. The Hardy space $H^1(\LU)$ is defined in terms of the maximal function for the semigroup $\K{t}{U} = \exp(-t\LU)$, namely
    $$H^1(\LU) = \set{f\in \LLL \ : \ \norm{f}_{H^1(\LU)}:= \norm{\sup_{t>0} \abs{\K{t}{U}f}}_\LLL}<\8.$$
Assume that $U=V+W$, where $V\geq 0$ satisfies the global Kato condition
$$\sup_{x\in \Rd} \int_{\Rd} V(y)|x-y|^{2-d} <\8.$$
We prove that, under certain assumptions on $W\geq 0$, the space $H^1(\LU)$ admits an atomic decomposition of local type. An atom $a$ for $H^1(\LU)$ is either of the form $a(x)=|Q|^{-1}\chi_Q(x)$, where $Q$ are special cubes determined by $W$, or $a$ satisfies the cancellation condition $\int_\Rd a(x)\w(x)\, dx = 0$, where $\w$ is an $(-\Deltax+V)$-harmonic function given by $\w(x) = \lim_{t\to \8} \K{t}{V}\mathbf{1}(x)$. Furthermore, we show that, in~some cases, the cancellation condition $\int_\Rd a(x)\w(x)\, dx = 0$ can be replaced by the classical one $\int_\Rd a(x)\, dx = 0$. However, we construct another example, such that the atomic spaces with these two cancellation conditions are not equivalent as Banach spaces.
\end{abstract}
\maketitle

\section{Background and statement of results}\label{intro}

\subsection{Introduction}
Let $U$ be a non-negative, locally integrable function on $\Rd$. In this article we consider the Schr\"odinger operator given by
    $
    -\Deltax + U,
    $
where $\Deltax$ is the standard Laplacian on $\Rd$ and $U$ is called {\it the potential}. Throughout the whole paper we assume that $d\geq 3$.

To be more precise, let us recall what do we mean by the Schr\"odinger operator. First, define a~quadratic form
    \eqx{
    \Qq^U(f,g) = \int_\Rd \nabla f(x)\overline{\nabla g(x)}\, dx + \int_\Rd U(x)f(x)\overline{g(x)}\, dx
    }
with the domain $\mathrm{Dom}(\Qq^U) = \set{f\in L^2(\Rd) \ : \ \nabla f, \sqrt{U} f \in L^2(\Rd)}.$ This quadratic form is closed, thus it defines the self-adjoin operator $\LU:\mathrm{Dom}(\LU) \to L^2(\Rd)$. In particular,
    \eqx{
    \mathrm{Dom}(\LU) =\set{f\in \mathrm{Dom}(\Qq^U) \ : \ \exists h\in L^2(\Rd) \ \forall g\in \mathrm{Dom}(\Qq^U) \quad  \Qq^U(f,g) = \int_\Rd h(x) \overline{g(x)}\, dx}
    }
and $\LU f := h$, when $f$ and $h$ are as above. Formally, we write $$\LU = -\Deltax + U.$$
Let $\expr{\K{t}{U}{}}_{t>0}$ be the semigroup generated by $\LU$ on $L^2(\Rd)$. By the Feynman-Kac formula, $\K{t}{U}$ has an integral kernel $\k{t}{U}{x,y}$ satisfying upper-Gaussian bounds, i.e.
    \eq{\label{up_gauss}
    0\leq \k{t}{U}{x,y} \leq (4\pi t)^{-d/2} \exp\expr{-\frac{|x-y|^2}{4t}} = P_t(x-y).
    }

The Hardy space $H^{1}(\LU)$ associated with $\LU$ is defined as follows. Let
    \eqx{\Mm^U f (x) = \sup_{t>0} \abs{\K{t}{U} f(x)}}
be a maximal operator associated with $\(\K{t}{U}\)_{t>0}$. We say that a function $f\in \LLL$ belongs to the maximal Hardy space $H^1(\LU)$, when
        \eq{\label{max}
        \norm{f}_{H^1(\LU)} := \norm{\Mm^U f (x)}_\LLL  <\8.
        }

In the paper atomic Hardy spaces play a special role. The general definition is as follows. Assume that a family of functions $\aa \subseteq \LLL$ is given. A function $a\in \aa$ will be called {\it an atom} and we assume that $\norm{a}_{\LLL} \leq 1$. We say that a function $f$ belongs to the atomic Hardy space $H^1_{at}(\aa)$, if
    \eq{\label{at_def}
    f(x) = \sum_{j=1}^\8 \la_j a_j(x),
    }
where  $a_j \in \aa, \ \la_j \in \CC, \ \text{and} \ \sum_{j=1}^\8 \abs{\la_j} <\8$.
Whenever $f \in H^1_{at}(\aa)$ we set
    \eq{\label{at_def_norm}
    \norm{f}_{H^1_{at}(\aa)}= \inf \set{\sum_{j=1}^\8 \abs{\la_j} \ : \ f \text{ as in (\ref{at_def})}}.
    }
It is not difficult to check that $H^1_{at}(\aa)$ is a Banach space and $H^1_{at}(\aa)\subseteq \LLL$.

In the classical theory of Hardy spaces an important result is the atomic decomposition theorem, see \cite{Coifman_Studia}, \cite{Latter_Studia}. It asserts that $H^1(-\Deltax)=H^1_{at}(\aa_{class})$ and the corresponding norms are equivalent. Here $\aa_{class}$ is the set of classical atoms, that is $a\in \aa_{class}$ if there exist a cube $Q$, such that $\supp\, a \subseteq Q$ (localization condition), $\norm{a}_\8 \leq |Q|^{-1}$ (size condition), and $\int_Q a(x) \, dx=0$ (cancellation condition). By $|S|$ we denote the Lebesgue measure of a set $S$ and
\eqx{
Q = Q(c_Q,r_Q) = \set{y=(y_1,...,y_d)\in \Rd \ : \ \max_{i=1,...,d}(|(c_Q)_i-y_i|)<r_Q},
}
where $c_Q$ and $r_Q$ are the center and the radius of $Q$, respectively. Denote $d_Q = \mathrm{diam}(Q) = 2\sqrt{d}r_Q$.

The question we shall be concerned with is: whether $H^1(\LU)$ coincides with $H^1_{at}(\aa)$ for a potential $U$ and a family $\aa$? If so, are the norms \eqref{max} and \eqref{at_def_norm} comparable?

There are partial answers to the question above. A general result of Hofmann et. al. \cite{Hofmann_Memoirs} gives an atomic and molecular characterizations of $H^1(\LU)$ for any positive potential $U \in L^1_{loc}(\Rd)$. Also, using \cite{Hofmann_Memoirs}, Dziubański and Zienkiewicz in \cite{DZ_Revista2} proved another general atomic characterization of $H^1(\LU)$. The atoms in \cite{Hofmann_Memoirs} are of the form $a=(\LU)^M b$, where $M\geq 1$ is fixed natural number and $b$ satisfies some localization and size conditions, see \cite[Theorem 7.1]{Hofmann_Memoirs}. Likewise, atoms in \cite{DZ_Revista2} are given by $a = \K{t}{U} b - b$ for similar $b$.

Although the approaches just mentioned are useful in many situations, they have also some disadvantages. One of them is that the atoms are images of the operator $\LU$ (or its semigroup) of some function, and they no more satisfy simple geometric conditions (localization, size, cancellation). One would also like to better understand the nature of $H^1(\LU)$ by describing it in terms of simpler, ''geometric atoms''. In the 90's Dziuba\'{n}ski and Zienkiewicz started studies on atomic decompositions of Hardy spaces for Schr\"odinger operators. In this paper we continue this approach. For more results of this type see \cite{CzajaZienkiewicz_ProcAMS}, \cite{Dziubanski_JFAA}, \cite{DP_Arkiv}, \cite{DP_Potential}, \cite{DP_JMAA}, \cite{DZ_Studia} \cite{DZ_Revista1}, \cite{DZ_Annali}, \cite{DZ_Revista2}, \cite{DZ_JFAA}. Let us finally mention, that this approach was successfully used e.g. for proving Riesz transform characterization of $H^1(\LU)$, while such characterization is not known in general.

Before proceeding to our main results, we present results of \cite{DZ_JFAA} and \cite{DZ_Studia}, which are the starting point for our considerations.

\subsection{The space $H^1(\LV)$} \label{subsecLV}
Assume that a potential $V\geq 0$ satisfies
    \al{\tag*{(S)}\label{S}
    \sup_{x\in \Rd} \int_\Rd  |x-y|^{2-d} V(y)  \, dy <\8.
    }
In other words, $(\LV)^{-1} V \in L^\8 (\Rd)$. Let $\w =\w(V)$ be a function defined by
    \eq{\label{omega}
    \omega(x) = \lim_{t\to \8}\int_\Rd \k{t}{V}{x,y} \, dy.
    }
The function $\w$ is $\LV$-harmonic and satisfies
    \eq{\label{omega_bounds}
    0<\delta<\omega(x) \leq 1,
     }
     with some $\delta$ for all $x\in \Rd$, see \cite[Lemma 2.1]{DZ_JFAA}. It is well-known, see \cite{Semenov}, that the integral kernel $\k{t}{V}{x,y}$ has not only upper-Gaussian bounds, but also lower-Gaussian bounds, that is we have $\kappa_1,\kappa_2>0$ such that
    \eq{\label{low_gauss}
    \k{t}{V}{x,y} \geq \kappa_1 t^{-d/2} \exp\(-\frac{\abs{x-y}^2}{\kappa_2 t}\)  .
    }

By definition, a function $a$ is an $\w$-atom, if there exists a cube $Q$ such that
    \eqx{
    \supp \, a \subseteq Q, \quad \|a\|_{\8} \leq |Q|^{-1}, \quad \text{and} \quad \int_Q a(x) \omega(x) \, dx =0.
    }
    Let $\aa_\w$ be the set of $\w$-atoms. Corollary 1.2 of \cite{DZ_JFAA} states that $H^1(\LV)= H^1_{at}(\aa_\w)$ and
    \eq{\label{DZ-JFAA}
     \| f\|_{H^1(\LV)} \simeq  \| f\|_{H^1_{at}(\aa_\w)}
    }

Let us mention that \ref{S} is satisfied for example when $V$ is compactly supported and $V \in L^p(\Rd)$ for some $p>d/2$. For more general examples, see \cite{DP_JMAA} and \cite{DZ_JFAA}.

\subsection{The space $H^1(\LW)$} \label{subsecLW}
For $\theta>0$ (small) and a cube $Q=Q(c_Q,r_Q)$ denote $Q^{*} = Q(c_Q, (1+\theta) r_Q)$. Assume a family of cubes $\qq$ is given and there exist $C,\theta>0$ such that for $Q_1, Q_2 \in \qq$, $Q_1 \neq Q_2$, we have:
    \al{\tag*{($G_1$)}\label{G1}
        &\textstyle{{\bigcup_{Q\in \qq} \mathrm{cl}(Q)} = \Rd},&&\\
        \tag*{($G_2$)}\label{G2}
        &|Q_{1} \cap Q_{2}| = 0,&&\\
        &\text{if }Q_{1}^{****}\cap Q_{2}^{****} \neq \emptyset  \text{, then } C^{-1} d_{Q_{1}} \leq d_{Q_{2}} \leq C d_{Q_{1}}.&&
        \tag*{($G_3$)}\label{G3}
    }
Observe that, under these assumptions, the family $\set{Q^{****} : Q\in \qq}$ is automatically a finite covering of $\Rd$. In the following, we shortly write that $\qq$ satisfies $\label{G}(G)$, when it satisfies \ref{G1}, \ref{G2}, \ref{G3}.

Suppose that for a potential $W\geq 0$ and a family $\qq$ as above there exist positive constants $\e, \delta, C$ such that
   \al{\tag*{(D)}\label{D}
    &\sup_{y\in Q^{**}} \int_\Rd \k{2^n d_Q^2}{W}{x,y}\, dx \leq Cn^{-1-\e} && (Q\in \qq, n \in \NN),\\
    \tag*{(K)}\label{K}
    &\int_0^{2t} (\mathbf {1}_{Q^{***}}W)*P_s(x)\, ds \leq C\left(\frac{t}{d_Q^2}\right)^\delta && (x\in \Rd,   \  Q\in \qq,t\leq d_Q^2),
    }
where $P_t(x-y) = \k{t}{0}{x,y}$ is the classical heat semigroup, see \eqref{up_gauss}. By definition, an $\qq$-atom is a function $a$ such that one of the following holds:
    \ite{
    \item $\text{there exists }  Q\in \qq \text{ and a cube } K\subset Q^{**}  \text{ such that: }$
    \eqx{
    \supp \, a \subseteq K, \quad \|a\|_{\8} \leq |K|^{-1}, \quad \int_{K} a(x)\, dx = 0;
    }
    \item $a(x) = |Q|^{-1} \chi_Q(x)$ for some $Q\in \qq$.
    }

    Let $\aa_\qq$ be a set of $\qq$-atoms. By Theorem 2.2 of \cite{DZ_Studia} we have that $H^1(\LW) = H^1_{at}(\aa_\qq)$ and
    \eqx{
     \| f\|_{H^1(\LW)} \simeq \| f\|_{H^1_{at}(\aa_\qq)}.
    }

    A list of examples of potentials $W$ and related families $\qq$ can be found in \cite{DZ_Studia}. At this place we shall only mention one simple example, that we shall use later in this paper. Let $t>0$ and denote by $\qq^{[t]}$ the family of cubes of radius equal to $t$ that satisfies $(G)$. If $W^{[t]}(x) = t^{-2}$, then the pair $(W^{[t]}, \qq^{[t]})$ satisfies \ref{D}, \ref{K}, $(G)$ with constants independent of $t$.

\subsection{Main results}
In this paper $V$ always denote a potential satisfying \ref{S} and $\w$ is related to $V$ by \eqref{omega}. Similarly, the pair $W$, $\qq$ always satisfy \ref{D}, \ref{K}, and $(G)$. Notice, that in $H^1_{at}(\aa_\w)$ and $H^1_{at}(\aa_\qq)$ two different effects appear. For an atom $a\in \aa_\w$ (atom for $\LV$) the cancellation condition is w.r.t. the measure $\w$, not the Lebesgue measure. On the other hand, for $a\in \aa_\qq$, there are ''local'' atoms, i.e. atoms of the type $|Q|^{-1} \chi_Q(x)$ that do not satisfy any cancellation condition.

The goal of this paper is to study $\LVW$ and its Hardy space $H^1(\LVW)$. We shall prove that in atomic decompositions for this space both effects described above appear simultaneously. Define $\aa_{\w,\qq}$ to be the set of $(\w,\qq)$-atoms, that is functions such that one of the following holds:
    \ite{
    \item $\text{there exists }  Q\in \qq \text{ and a cube } K\subset Q^{**}  \text{ such that: }$
    \eqx{
    \supp \, a \subseteq K, \quad \|a\|_{\8} \leq |K|^{-1}, \quad \int_{K} a(x)\w(x)\, dx = 0,
    }
    \item $a(x) = |Q|^{-1} \chi_Q(x)$ for some $Q\in \qq$.
    }

The following theorem gives the atomic characterization of $H^1(\LVW)$ in the spirit of \cite{DZ_Studia} and \cite{DZ_JFAA}.
    \sthm{mainthm}{
        Assume that $d\geq 3$, $V\geq 0$ satisfies \ref{S}, and $W\geq 0$ with a family $\qq$ satisfy \ref{D}, \ref{K}, $(G)$. Then
        \eq{\label{mainest}
        C^{-1} \| f\|_{H^1(\LVW)} \leq \| f\|_{H^1_{at}(\aa_{\w,\qq})} \leq C \| f\|_{H^1(\LVW)}.
        }
        In particular, $H^1(\LVW) = H^1_{at}(\aa_{\w,\qq})$.
    }

In Theorem \ref{mainthm} atoms are localized to cubes $Q\in \qq$ and the cancellation condition is w.r.t. the measure $\w(x)\,dx$. However, it is not hard to see that every $(\w,\qq)-atom$ can be written as a linear combination of just $\qq$-atoms. Indeed, if $a$ is such that $\supp a \subseteq K \subseteq Q^{**}$, $\norm{a}_\8 \leq |K|^{-1}$, and $\int_{K} a(x) \w(x)\, dx =0$ for $Q\in \qq$, then
    \eqx{
    a(x) = \expr{a(x)-\kappa |Q|^{-1} \1_{Q}(x)} + \kappa |Q|^{-1} \1_{Q}(x) = b_1(x)+b_2(x),
    }
where $\kappa=\int_{K}a(x)\, dx$, $|\kappa| \leq 1$. Observe that $\supp \, b_1 \subseteq  Q^{**}$ and $\int_{Q^{**}} b_1(x)\, dx=0$. Thus both $b_1$ and $b_2$ are multiples of $\qq$-atoms. What we have just shown is that $a\in H^1_{at}(\aa_\qq) $ and
    \eq{\label{T}
    \norm{a}_{H^1_{at}(\aa_\qq)}\leq T,
    }
for every $(\w,\qq)$-atom $a$.

The constant $T$ in \eqref{T} possibly depend on $a$. This lead us to the following question: whether $H^1_{at}(\aa_{\w,\qq})$ and $H^1_{at}(\aa_{\qq})$ are equal as Banach spaces? In Theorem \ref{thm2} we prove that, under certain Lipshitz assumption, the answer to this question is positive. However, a more difficult task is to find an example such that $\norm{f}_{\HQ} \not\simeq \norm{f}_{\HwQ}$. This is done in Example \ref{C}.

    \sthm{thm2}{
    Assume that $0<\delta\leq \w\leq 1$, $\qq$ satisfies $(G)$, and there exists $\lambda >0$ such that
        \eq{\label{lipshitz}
        |\w(x) - \w(y)| \leq C \expr{ \frac{|x-y|}{d_Q}}^\lambda \qquad (Q\in \qq, x,y\in Q^{**}).
        }
    Then
        \eq{\label{comp}
        \norm{f}_{\HQ} \simeq \norm{f}_{\HwQ}.
        }
    }

As an example that fulfills the assumptions of Theorem \ref{thm2} one could take $W^{[1]}$, $\qq^{[1]}$  (see Subsection \ref{subsecLW}) and $\w = \w(V)$, with $V$ such that $\supp \, V \subseteq Q(0,1)$ and $V \in L^p(\Rd)$ for $p>d/2$ (for details see \cite{DZ_Annali}). In this case $\w$ satisfies global H\"older condition.

    \sex{example}{\label{C}
     Let $\qq^{[1]}$ be as above, and $\w=\w(\vv)$, where $\vv$ is a potential given in \eqref{V}. There exist a sequence of $(\w,\qq)$-atoms $a_j$, such that
        \eq{\lab{not-comparable}
        \limj \norm{a_j}_{\HQ} = \8.
        }
    In other words, $\norm{f}_{\HwQ} \not \simeq \norm{f}_{\HQ}$.
    }

The paper is organized as follows. Section \ref{sec2} is devoted to local Hardy spaces. We prove an atomic decomposition for a local version of $H^1(\LV)$. In Section \ref{sec3} we prove some auxiliary estimates, most of which are analogues of Lemmas in \cite{DZ_Studia}. In Section \ref{sec4} and Section \ref{sec5} we present the proofs of Theorems \ref{mainthm} and \ref{thm2}, respectively. In Section \ref{sec6} we provide details of Example \ref{example} and prove \eqref{not-comparable}. Finally, in the Appendix we give a proof of $\norm{f}_\LLL \leq \norm{\sup_{t\leq \tau} \K{t}{U}f}_\LLL$.

At the end of this section let us give a short remark. In some papers authors define local atomic spaces in a slightly different manner. The remark below clarify, that different definitions lead to the same atomic Hardy spaces in the sense of equivalent Banach spaces.

\rem{remX}{
Let us consider $\qq$ and $\w$ as above and a function $\mathfrak a$ that satisfies:
    \ite{
    \item $\text{there exists }  Q\in \qq \text{ and a cube } K\subset Q^{**}  \text{ such that: }$
    \eqx{
    \supp \, \mathfrak a \subseteq K, \quad 4d_K\geq d_Q, \quad  \|\mathfrak a\|_{\8} \leq |K|^{-1}.
    }
    }
 For each $\mathfrak a$ as above, we have that $\norm{\mathfrak a}_{H^1_{at}(\aa_{\qq,\w})} \leq C$, with universal $C$. To see this, one has to write $\mathfrak a$ as a linear combination of $|Q|^{-1}\chi_Q(x)$ and atom with cancellation condition. Therefore, the functions $\mathfrak a$ as above can be substitutes for the atoms of the form $|Q|^{-1}\chi_Q(x)$ in the definition of $\aa_{\w,\qq}$.
}

\section{Local Hardy spaces}\label{sec2}
\subsection{Local Hardy spaces}

In this section we put aside $W$ and $\qq$ for a moment and consider only $\LV$ and related objects. The local version of the maximal operator $\Mm^V$ at scale $\tau >0$ is
        \eqx{\label{max-loc}
        \Mm^{V}_{\tau} f (x) = \sup_{t\leq\tau^2} \abs{\K{t}{V} f (x)}.
        }
By definition, a function $f\in \LLL$ is in the local Hardy space $h_\tau^1(\LV)$, when $\Mm^{V}_{\tau} f$ is in $\LLL$. We set
    \eqx{
    \|f\|_{h_\tau^1(\LV)} := \|\Mm^{V}_{\tau} f\|_{\LLL}.
    }

In a special case $V \equiv 0$, the space $h^1_\tau(-\Deltax)$ is a classical local Hardy space introduced by Goldberg \cite{Goldberg_Duke}. It follows from \cite{Goldberg_Duke} that
   \eq{\label{Gold}
    C^{-1} \| f\|_{H^1_{at}(\aa_{\qq^{[\tau]}})}\leq \| f\|_{h_\tau^1(-\Deltax)} \leq C \| f\|_{H^1_{at}(\aa_{\qq^{[\tau]}})},
    }
where $C$ does not depend on $\tau$. The following proposition is a generalization of \eqref{Gold} for $h^1_\tau (\LV)$ localized to a cube of diameter comparable to $\tau$. It will play a crucial role in the proof of Theorem~\ref{mainthm}.
    \prop{lem_local}{Let $Q$ be a cube.
    \en{[a)]
        \item
            Let $a$ be $\w-$atom, such that $\supp \, a \subseteq Q^{**}$ {\bf or} $a(x)= |Q|^{-1} \chi_Q(x)$. Then
                \eq{\lab{on_atom_loc}
                \norm{\Mm^V_{d_Q} a}_{\LLL} \leq C.
                }
        \item
            Assume that $\supp \, f \subseteq Q^*$and
            $
            \Mm^{V}_{d_Q} f(x)
            \in L^1(\Rd ).
            $
            There exist $\la_j$ and $a_j$ being either $\w$-atoms {\bf or} of the form $ |Q|^{-1} \chi_Q(x)$, such that
                \eqx{
                f(x) = \sumj \la_j a_j(x), \quad \quad \sum_{j=1}^\8 \abs{\la_j} \leq C \norm{\Mm_{d_Q}^{V}f}_{L^1\(\Rd\)} .
                }
            The constant $C$ above depends only on $d$ and $\theta$ in the definition of $Q^*$.
    }
    }

    \pr{
    Assume first that $a$ is $\w-$atom. Obviously, $\Mm^V_{d_Q} a(x) \leq \Mm^V a(x)$, so \eqref{on_atom_loc} holds by \eqref{DZ-JFAA}. In the case when $a = |Q|^{-1} \chi_Q$ we use \eqref{up_gauss} and \eqref{Gold}, getting
        \eqx{
        \norm{\Mm^V_{d_Q} a}_\LLL \leq \norm{\Mm^0_{d_Q} a}_\LLL \leq C.
        }

    Now, let $f$ be as in the assumptions of b). Set
        \eqx{
        g(x) = f(x) - \K{d_Q^2/2}{V}f(x),
        }
    so that
        \eqx{
        f(x) \w(x) = g(x) \w(x) + \K{d_Q^2/2}{V}f(x) \w(x) =h_1(x) + h_2(x).
        }
    We claim that $h_1 \in H^1(-\Deltax)$ and $h_2 \in h^1_{d_Q}(-\Deltax)$ with
        \al{\label{est5}
        &\norm{h_1}_{H^1(-\Deltax)} \leq C \norm{\Mm^V_{d_Q}f}_{\LLL},\\
        \label{est6}
        &\norm{h_2}_{h^1_{d_Q}(-\Deltax)} \leq C \norm{\Mm^V_{d_Q}f}_{\LLL}.
        }
    To prove \eqref{est5}, observe that
        \eqx{
        \norm{\sup_{t\leq d_Q^2/2}\abs{\K{t}{V} g}}_\LLL \leq 
         2 \norm{\Mm^V_{d_Q}f}_{\LLL} <\8.
        }
    Likewise, $$\norm{\sup_{t> d_Q^2/2}\abs{\K{t}{V} g(x)}}_\LLL\leq C\norm{f}_\LLL$$ by the argument identical as in the proof of \cite[Proposition 6.3]{DZ_JFAA}. By Corollary \ref{cor_pointwise}, $\norm{f}_\LLL \leq \norm{\Mm_{d_Q}^V f}_\LLL$. Thus $g \in H^1(\LV)$ and, by \eqref{DZ-JFAA}, $h_1 =g\cdot \w\in H^1(-\Deltax)$, so \eqref{est5} is proved.

    Now, we turn to prove \eqref{est6}. It is clear that
        \eqx{
        h_2(x) = \sum_{K\in \qq^{[d_Q]}} \K{d_Q^2/2}{V}f(x) \w(x) \chi_K(x) = \sum_{K\in \qq^{[d_Q]}} h_K(x)
        }
    and
    \spx{
    \norm{h_K}_\8 &\leq C \int_K d_Q^{-d} \exp\expr{-\frac{|x-y|^2}{2d_Q^2}} \abs{f(y)}\, dy\\
        &\leq C |Q|^{-1} \exp\expr{-\frac{d(Q^*,K)^2}{2d_Q^2}}  \norm{\Mm^V_{d_Q}f}_{\LLL}.
    }
    Clearly, $\supp \, h_K\subseteq K$, so by using the classical atomic characterization of $h^1_{d_Q} (-\Deltax)$ we have that $\norm{h_K}_{h^1_{d_Q}(-\Deltax)} \leq C \exp\expr{-\frac{d(Q^*,K)^2}{2d_Q^2}}  \norm{\Mm^V_{d_Q}f}_{\LLL}$. Summing up,
        \spx{
        \norm{h_2}_{h^1_{d_Q}(-\Deltax)} 
         &\leq C \norm{\Mm^V_{d_Q}f}_{\LLL} \sum_{K\in \qq^{[d_Q]}} \exp\expr{-\frac{d(Q^*,K)^2}{2d_Q^2}}\\
         &\leq C \norm{\Mm^V_{d_Q}f}_{\LLL},
        }
    where the last inequality is a simple geometric observation.

    Having \eqref{est5} and \eqref{est6} proved, we finish the prove by the following argument. The function $f\cdot \w$ is supported in $Q^*$ and $f\cdot \w \in h^1_{d_Q} (-\Deltax)$ with $\norm{f\cdot \w}_{h^1_{d_Q}(-\Deltax)}\leq C \norm{\Mm^V_{d_Q}f}_{\LLL}$. So, by the classical local characterization of $h^1_{d_Q}(-\Deltax)$, $f\cdot \w = \sum_j \la_j a_j$, where $a_j$ are either classical atoms or of the form $|Q|^{-1} \chi_Q(x)$. Moreover, $\sum_j \abs{\la_j} \leq C \norm{\Mm^V_{d_Q}f}_{\LLL}$. Then $f = \sum_j \la_j b_j$ where $b_j = a_j/\w$ are either $\w-$atoms or $b_j = \w^{-1} |Q|^{-1} \chi_Q$. In the last case, $b_j$ can be decomposed into a linear combination of $|Q|^{-1}\chi_Q$-atom and $\w$-atom, exactly as in Remark~\ref{remX}.
    }

The following corollary is a ''global'' version of Proposition \ref{lem_local} and can be proved by standard techniques. The details are left to the reader.

    \cor{prop_local}{
    There exists a constant $C$, independent of $\tau>0$, such that
        \eqx{
         \| f\|_{h_\tau^1(\LV)} \simeq \| f\|_{H^1_{at}(\aa_{\w,\qq^{[\tau]}})}.
        }
    In particular, $h_\tau^1(\LV) = H^1_{at}(\aa_{\w,\qq^\tau})$.
    }

\section{Auxiliary estimates}\label{sec3}

In this section we present tools and lemmas that will be used in the proof of Theorem~\ref{mainthm}. The proofs of Lemmas \ref{global_part}, \ref{perturbation}, \ref{commutator}, \ref{another_global} are very similar to their analogues in \cite{DZ_Studia}. Thus we only provide sketches how to adapt proofs from \cite{DZ_Studia} to our background.

Let $U_1,U_2\geq 0$ be two potentials. A well-known perturbation formula states that
    \eq{\label{pert}
    \K{t}{U_1} - \K{t}{U_1+U_2} = \int_0^t \K{t-s}{U_1} \, U_2 \, \K{s}{U_1+U_2} \, ds.
    }
For the kernels this reads as
    \eq{
    \label{pert_ker}
    \k{t}{U_1}{x,y} - \k{t}{U_1+U_2}{x,y} = \int_0^t \int_{\Rd} \k{t-s}{U_1}{x,z}V(z) \k{s}{U_1+U_2}{z,y}\, dz \, ds.
    }

With a family $\qq$ satisfying $(G)$ we associate a partition of unity $\Phi = \{\phi_Q\}_{Q\in \qq}$ such that
    \eq{\label{partition}
    0\leq \phi_Q \in C_c^{\8} ( Q^{*}), \quad {\bf 1}_{\Rd} = \sum_{Q\in \qq} \phi_Q, \quad \norm{\nabla \phi_Q}_\8 \leq C d_Q^{-1}.
    }

    \lem{global_part}{
    Let $U\in L^1_{loc}(\Rd)$ be a positive potential. For $f\in \LLL$ and $Q \in \qq$,
        \eqx{
        \norm{ \sup_{t\leq d_Q^2} \abs{\K{t}{U}(\phi_Qf)}}_{L^1((Q^{**})^c)} \leq \norm{\phi_Q f}_{\LLL}.
        }
    }

    \pr{
    Let $c_Q$ be the center of $Q$. For $t\leq d_Q^2$, $y\in Q^{*}$ and $x \not \in Q^{**}$ we have
        \alx{
        \sup_{t\leq d_Q^2} \k{t}{U}{x,y} \leq \sup_{t\leq  d_Q^2} C t^{-d/2} \exp\(-\frac{|x-c_Q|^2}{c t}\) \leq C d_Q^{-d} \exp\(-\frac{|x-c_Q|^2}{c d_Q^2}\).
        }
     The lemma follows by integrating the last expression w.r.t. $dx$ on $(Q^{**})^c$.
    }

    \lem{perturbation}{
    Assume \ref{K}. For $f\in \LLL$ and $Q \in \qq$,
        \eqx{
        \norm{ \sup_{t\leq d_Q^2} \abs{\(\K{t}{V} - \K{t}{V+W}\) (\phi_Q f)} }_{\LLL}\leq  C \norm{\phi_Q f}_\LLL.
        }
    }

    \pr{[Sketch of the proof] Using \eqref{pert} we write
        \spx{
        \(\K{t}{V} - \K{t}{V+W}\) (\phi_Q f) &=\int_0^t \K{t-s}{V} (W \cdot \1_{(Q^{***})^c}) \K{s}{V+W}\(\phi_Q f\) \,ds\cr
        &\quad + \int_0^t \K{t-s}{V} (W \cdot \1_{Q^{***}}) \K{s}{V+W}\(\phi_Q f\) \,ds.
        }
    Both summands can be estimated similarly as in \cite[Lemma 3.11]{DZ_Studia}. In order to repeat arguments of \cite{DZ_Studia}, one should have in mind that, by \eqref{pert_ker},
    \eq{\label{monotonicity}
    \k{t}{V+W}{x,y} \leq \k{t}{U}{x,y}\leq P_t(x-y),}
    where $U$ is either $V$ or $W$. The details are omitted.
    }

For each $Q\in \qq$ we set
    \alx{
        &\qq_{loc, Q} = \{ Q' \in \qq \ : \ Q^{***} \cap Q'^{***} \neq \emptyset\}, \\
        &\qq_{glob,Q} = \{ Q'' \in \qq \ : \ Q^{***} \cap Q''^{***}  =\emptyset\}.
    }
Roughly speaking, for each $Q$, the set $\qq_{loc, Q}$ is the set of cubes $Q' \in \qq$ that are ''close'' to $Q$. For a function $f$ denote
    \eqx{
    f_{loc,Q}= \sum_{Q'\in \qq_{loc,Q}} \phi_{Q'} f, \quad f_{glob,Q} = f - f_{loc,{Q}}.
    }
The following two lemmas and their proofs are almost identical to \cite[Lemma 3.7]{DZ_Studia} and \cite[Lemma 3.8]{DZ_Studia}. To see this one only has to use \eqref{monotonicity}. The details are left to the reader.
    \lem{commutator}{
    For $f\in \LLL$ and $Q\in \qq$,
        \eqx{
        \norm{ \sup_{t>0} \abs{ \K{t}{V+W}(\phi_Q \cdot f_{loc,Q}) - \phi_Q \cdot \K{t}{V+W}(f_{loc,Q}) }}_{L^1(Q^{**})} \leq C \norm{f_{loc,Q}}_{\LLL}.  
        }
    }

    \lem{another_global}{
    Assume \ref{D}. For $f\in \LLL$ and $Q\in \qq$,
        \eqx{
        \sum_{Q\in \qq} \norm{ \sup_{t\leq d_Q^2}\abs{ K_t^{V+W} (f_{glob, Q})}}_{L^1(Q^{*})} \leq C\norm{f}_\LLL.
        }
    }

\section{Proof of Theorem \ref{mainthm}}\label{sec4}

In the proof below, we shall often use the fact that, for $0\leq U \in L^1_{loc}(\Rd)$ and $\tau>0$, we have
    \eq{
    \label{app}
    \norm{f}_{\LLL} \leq \norm{\Mm_\tau^U f}_{\LLL}.
    }
This is a consequence of semigroup property and Gaussian estimates. A detailed proof is given in the Appendix, see Proposition \ref{prop_pointwise} and Corollary \ref{cor_pointwise}.

{\bf First implication.} We start by proving the second inequality of \eqref{mainest}, that is for a function $f$ such that $\norm{f}_{H^1(\LVW)} <\8$ we will find $(\w,\qq)$-atoms $a_i$ such that
    \eqx{
    f(x) = \sumi \la_i a_i(x) \quad \text{and} \quad \sumi \abs{\la_i} \leq C \norm{f}_{H^1(\LVW)}.
    }
Let $\phi_Q$ be as in \eqref{partition}, in particular $f = \sum_{Q\in \qq} \phi_Qf$. The key estimate is the following.
    \eq{\label{key_first_imp}
    \sum_{Q\in \qq} \norm{\sup_{t\leq d_{Q}^2} \abs{\K{t}{V} ( \phi_Q f) (x)}}_{L^1\({\Rd}\)} \leq C \norm{f}_{H^1(\LVW)}.
    }

Now we prove \eqref{key_first_imp}. By Lemma \ref{global_part} we get that $\sum_{Q\in \qq} \norm{\cdot}_{L^1((Q^{**})^c)} \leq C \norm{f}_{\LLL}$. Now we concentrate our attention on $Q^{**}$. Notice that
    \alx{
    \K{t}{V} (\phi_Qf) = &\left[\(\K{t}{V} - \K{t}{V+W}\) ( \phi_Q f) \right] + \left[\K{t}{V+W}(\phi_Q f) - \phi_Q \cdot \K{t}{V+W}(f_{loc,Q})  \right] \\
    &+ \left[ - \phi_Q \cdot \K{t}{V+W}(f_{glob,Q})\right]+\left[\phi_Q \cdot \K{t}{V+W}(f)\right]\\
    =& A_1+A_2+A_3+A_4.
    }
Notice that $\phi_Q f_{loc,Q} = \phi_Q f$. Lemmas \ref{perturbation}, \ref{commutator}, \ref{another_global} lead to
    \alx{
    \sum_{k=1}^3 \sum_{Q\in \qq} \norm{ \sup_{t\leq d_Q^2} \abs{A_k}}_{L^1(Q^{**})} & \leq C \sum_{Q\in\qq} \( \norm{\phi_Q\cdot f}_\LLL + \norm{f_{loc,Q}}_\LLL \) + \norm{f}_\LLL\\
    &\leq C \norm{f}_{\LLL} \leq C \norm{f}_{H^1(\LVW)},
    }
where we have used \eqref{app} and
    \spx{
    \sum_{Q\in\qq} \norm{f_{loc,Q}}_\LLL &\leq \sum_{Q\in\qq} \sum_{Q'\in\qq_{loc,Q}}\norm{\phi_{Q'}f}_{\LLL} = \sum_{Q'\in\qq} \sum_{Q\in\qq_{loc,Q'}}\norm{\phi_{Q'}f}_{\LLL}\\
    &\leq C \sum_{Q'\in\qq} \norm{\phi_{Q'}f}_{\LLL} \leq C \norm{f}_\LLL.
    }
The proof of \eqref{key_first_imp} is finished by noticing that
    \eqx{
    \sum_{Q\in \qq} \norm{ \sup_{t\leq d_Q^2} \abs{A_4}}_{L^1(Q^{**})} \leq C \norm{f}_{H^1(\LVW)}.
    }

Having \eqref{key_first_imp} proved, we apply Proposition \ref{lem_local}b to $\phi_Q f$, obtaining $\lambda_{j,Q}$ and $(\w,\qq)-$atoms $a_{j,Q}$ such that
    \eqx{
    \phi_Q(x) f(x) = \sum_{j=1}^\8 \la_{j,Q} a_{j,Q}(x), \quad \text{with} \quad \sum_{j=1}^\8 \abs{\la_{j,Q}} \leq C \norm{\sup_{t\leq d_{Q}^2} \abs{\K{t}{V} ( \phi_Q f) (x)}}_{\LLL}.
    }
Therefore,
    \eqx{
    f(x) = \sum_{j,Q} \la_{j,Q} a_{j,Q}(x), \ \ \text{with} \ \  \sum_{j,Q} \abs{\la_{j,Q}} \leq C \norm{f}_{H^1(\LVW)}
    }
and the proof of the first part is finished.

{\bf Second implication.} By a standard argument it is enough to prove that
    \eqx{
    \norm{\sup_{t>0} \abs{\K{t}{V+W}{a}}}_{\LLL} \leq C
    }
    for $a\in \aa_{\w,\qq}$. Assume then that $\supp \, a \subseteq Q^{**}$, where $Q \in \qq$. By the definition of $\qq_{loc,Q}$ and $\phi_Q$ it is clear that $a = a_{loc, Q}$. From \ref{G3} there exists a universal constant $m \in \NN$ such that $d_{Q'}^2 \geq 2^{-m} d_Q^2$ whenever $Q' \in \qq_{loc,Q}$.
    \spx{
    \norm{\sup_{t\leq 2^{-m}d_Q^2} \abs{\K{t}{V+W}{a}}}_{\LLL} \leq &\sum_{Q' \in \qq_{loc,Q}} \norm{\sup_{t\leq d_{Q'}^2} \abs{\(\K{t}{V+W} - \K{t}{V}\) (\phi_{Q'} a)}}_{\LLL} \\
    &+ \norm{\sup_{t\leq d_Q^2} \abs{\K{t}{V}{a}}}_{\LLL}
    }
By Lemma \ref{perturbation}, the sum is bounded by $C\norm{a}_{\LLL} \leq C$. The second summand is bounded by Proposition \ref{lem_local}a.

What is left is to consider $t\geq 2^{-m} d_Q^2$. Denote
    \eqx{I_j = [2^j d_Q^2, 2^{j+1} d_Q^2], \qquad I_j^\Diamond = [2^{j-1}d_Q^2, 3 \cdot 2^{j-1}d_Q^2].
    }
Note that $I_j = \set{x+2^{j-1}d_Q^2 \ : \ x\in I_j^\Diamond}$. By \eqref{up_gauss} it is not hard to check that for $g\in \LLL$ we have
    \eqx{
    \norm{\sup_{t\in I_j^\Diamond \cup I_j} \abs{\K{t}{V+W}g}}_\LLL \leq C \norm{g}_\LLL,
    }
where $C$ does not depend on $j$ and $g$. Therefore, for $j\geq 2$,
    \spx{
    \norm{\sup_{t\in I_j}\abs{\K{t}{V+W}a}}_{\LLL} &\leq \norm{\sup_{t\in I_j^\Diamond}\K{t}{V+W}\expr{\K{2^{j-1}d_Q^2}{W}\abs{a}}}_{\LLL}\\
    &\leq C\norm{\K{2^{j-1}d_Q^2}{W}\abs{a}}_\LLL \leq Cj^{-1-\e},
    }
where in the last inequality we have used \ref{D}. The proof is finished by noticing that
    \eqx{
    \norm{\sup_{t\geq 2^{-m}d_Q^2}\abs{\K{t}{V+W}a}}_{\LLL} \leq \sum_{j=-m}^\8 \norm{\sup_{t\in I_j}\abs{\K{t}{V+W}a}}_{\LLL}\leq C \expr{m+2+\sum_{j=2}^\8 j^{-1-\e}}\leq C.
    }

\section{Proof of Theorem \ref{thm2}}\label{sec5}

The proof follows by known procedure that uses atomic decompositions. Assume that $W,V,\qq, \w$ are given and $\w$ satisfies \eqref{lipshitz}.

To prove one of the inequalities of \eqref{comp} it is enough to show that
    \eq{\label{comparable}
    \norm{a}_{H^1_{at}(\aa_{\qq})}\leq C
    }
for $a\in \aa_{\w,\qq}$. Obviously, if $a$ is an atom of the form $a(x)=|Q|^{-1}\chi_Q(x)$, the inequality \eqref{comparable} holds with $C=1$. Assume then that $a$ is such that $\supp \, a \subseteq K \subseteq Q^{**}$, $Q\in \qq$, $\norm{a}_\8 \leq |K|^{-1}$, $ \int_{K}a(x)\w(x)\, dx=0$. Take a sequence of cubes $G_n$ such that
    \eqx{
    K = G_0 \subseteq G_1 \subseteq ... \subseteq G_N \subseteq Q^{**}, \quad d_{G_{n+1}} = 2 d_{G_n} \quad (n=0,...,N-1),
    }
and $d_Q\leq 2 d_{G_N}$. Observe that $N \leq C(\log_2 (d_Q/d_{K})+1)$ and $a(x) = \sum_{n=0}^{N+2} b_n(x)$, where
    \alx{
    &b_0(x) = a(x) - t_0 \chi_{G_0}(x),\\
    &b_n(x) = t_{n-1} \chi_{G_{n-1}}(x) - t_{n}\chi_{G_n}(x) \qquad (n=1,...,N),\\
    &b_{N+1} = t_N \chi_{G_N}(x) - t_{N+1} |Q|^{-1} \chi_Q(x),\\
    &b_{N+2} = t_{N+1} |Q|^{-1} \chi_Q(x).
    }
The constants $t_n$, are chosen so that $\int b_n(x) \, dx=0$ for $n=0,...,N+1$, i.e.
    \alx{
    &t_0 = |G_0|^{-1} \int_{G_0}a(x)\, dx,\\
    &t_{n} =  2^{-d}t_{n-1} \qquad (n=1,...,N),\\
    &t_{N+1} = t_N |G_N|.
    }
The key estimate, that uses \eqref{lipshitz} and the cancellation property, is the following
    \spx{
    \abs{t_0} &= \abs{K}^{-1} \w(c_{K})^{-1} \abs{\int_{K}a(x)\expr{\w(c_{K}) - \w(x)}\, dx}\\
    & \leq C \abs{K}^{-2} \int_{K}\expr{\frac{|x-c_K|}{d_Q}}^\la \, dx \leq C |K|^{-1}\expr{\frac{d_{K}}{d_Q}}^\la \leq C  2^{-c N} |K|^{-1}
    }
Thus $\abs{t_n} \leq C 2^{-c N} |G_n|^{-1}$ for $n=1,...,N$, and $\abs{t_{N+1}} \leq C$.

Obviously, $\supp \, b_n \subseteq G_n$ for $n=0,...,N$, and $\supp \, b_{N+1} \subseteq Q^{**}$. Moreover,
    \alx{\label{b_zero}
    &\norm{b_0}_\8 \leq \abs{K}^{-1} + \abs{t_0} \leq C \abs{K}^{-1},\\
    &\norm{b_n}_\8 \leq C |t_{n-1}| \leq  C 2^{-c N} \abs{G_n}^{-1} \quad (n=1,...,N)\\
    &\norm{b_{N+1}}_\8 \leq C|Q^{**}|^{-1}.
    }
As a consequence we have that all $b_n$ are multiples of $ H^1_{at}(\aa_{\qq})$-atoms and \eqref{comparable} is proved, since
    \eqx{
    \norm{a}_{H^1_{at}(\aa_{\qq})} \leq \sum_{n=0}^{N+2} \norm{b_n}_{H^1_{at}(\aa_{\qq})} \leq CN2^{-cN} +3C \leq C.
    }

For the second inequality one should consider $a\in H^1_{at}(\aa_\qq)$ and prove that
    \eqx{
    \norm{a}_{H^1_{at}(\aa_{\qq,\w})}\leq C.
    }
This can be done in a similar fashion. The details are omitted here.

\section{Example \ref{example}}\label{sec6}

Denote $c_n = 2^{n}\mathbf{e_1}$ and $C_n = Q(c_n, 1/(2n))$, where $\mathbf{e_1}$ denotes the vector $(1,0,...,0)$ in $\Rd$. The potential $\vv$ that we need for Example \ref{example} is the following

\eq{\label{V}
\vv(x) = \sum_{k=2}^\8 k^2 \chi_{C_k}(x).
}

\lem{l1}{
    $\vv$ satisfies \ref{S}.
}

\pr{
    Let $x\in \Rd$.
    \eqx{
    \int_\Rd \vv(y)|x-y|^{2-d} \, dy = \sum_{k=2}^\8 k^2 \int_{C_k} |x-y|^{2-d} \, dy= \sum_{k=2}^\8 I_k.
    }
    We have that
    \al{\label{est1}
    &I_k \leq k^2 \int_{C_k} |y-c_k|^{2-d} \, dy \leq C \qquad (x\in \Rd),\\
    \label{est2}
    &I_k \leq C k^2 \int_{C_k} |x-c_k|^{2-d} \, dy \leq C (k|x-c_k|)^{2-d} \qquad  (x\notin 2C_k).
    }
    Consider $x=(x_1, ... , x_d)$ and let $N\geq 2$ be such that $2^N< x_1 \leq 2^{N+1}$ ($N=2$ when $x_1 \leq 8$). Then
    \eqx{
    \sum_{k=2}^\8 I_k = \sum_{k=2}^{N-1} I_k+ \expr{I_{N} + I_{N+1}} +\sum_{k=N+2}^\8 I_k = A_1+A_2 +A_3,
    }
    with obvious modification when $N=2$. Obviously, $A_2\leq C$ by \eqref{est1}. Moreover, for $k\neq N$ and $k\neq N+1$, we have that $|x-c_k| \geq c2^{\max(N,k)}$, so using \eqref{est2} we obtain
    \alx{
    &A_1 \leq C \sum_{k=2}^{N-2} \expr{k2^N}^{2-d} \leq C,\\
    &A_3 \leq C \sum_{k=N+1}^\8 \expr{k2^k}^{2-d} \leq C.
    }
}

For the rest of this section by $\w$ we mean $\w(\vv)$ for $\vv$ given by \eqref{V}. The following lemma give an essential information about local oscillations of $\w$.

\prop{p1}{
    Let $c_n$ and $C_n$ be as above,
        \eqx{
        d_n = c_n + (\tau/n) \mathbf{e_1},  \quad D_n =Q(d_n, 1/(2n)).
        }
        There exists $\tau>3$, $c_0>0$, and $N\in \Nn$ such that for $n\geq N$ we have
        \eq{\label{example_osc}
        \inf_{x\in D_n, \ y\in C_n} \expr{\w(x) - \w(y)} \geq c_0.
        }
}

Let us remark that $\w$ satisfying \eqref{example_osc} cannot fulfill the global H\"older condition. To see this, just observe that $|c_n- d_n| \to 0$ and $\w(d_n) - \w(c_n) \geq c_0$.

\pr{
    Recall that $\k{t}{\vv}{x,y}$ always satisfies upper-Gaussian bounds, see \eqref{up_gauss}. By Lemma \ref{l1}, there are also lower-Gaussian bounds. Set $\kappa = \min(\kappa_1, \kappa_2)$, where $\kappa_1, \kappa_2$ are as in \eqref{low_gauss}. Put $U_1 = 0, \ U_2=\vv$ in \eqref{pert_ker}, integrate w.r.t. $x\in \Rd$, and let $t$ tend to infinity. We obtain that
        \eqx{
        1- \w(y) = \int_0^\8 \int_\Rd \vv(z) K_s^{\vv}(z,y)\, dy\, ds.
        }

    It is enough to show that, for properly chosen $\tau$ and $c_0$, the following estimates hold for $x\in D_n$ and $y\in C_n$.
    \eq{\label{est3}
        1-\w(y) =\int_0^\8 \int_\Rd \vv(z) K_s^{\vv}(z,y)\, dy\, ds \geq 2c_0,
    }
    \eq{\label{est4}
        1-\w(x) =\int_0^\8 \int_\Rd \vv(z) K_s^{\vv}(z,x)\, dy\, ds \leq c_0.
    }

    Fix $n\geq 2$ and $y\in C_n$. By \eqref{low_gauss} and \eqref{V},
        \alx{
            \int_0^\8 \int_\Rd \vv(z) \k{s}{\vv}{z,y}\, dy\, ds &\geq  c \int_0^\8 \int_{C_n} \kappa n^2 s^{-d/2} \exp\(-\frac{|z-y|^2}{\kappa s}\) \, dz \, ds\\
            &= c \kappa n^2 \int_{C_n} |z-y|^{2-d} dz \cdot \int_0^\8 s^{-d/2} \exp\(-\frac{1}{\kappa s}\) \, ds \\
            &\geq c(d,\kappa) =: 2c_0.
        }
    Thus \eqref{est3} is proved. For $x\in D_n$,
        \alx{
        \int_0^\8 \int_\Rd \vv(z) K_s^{\vv}(z,x)\, dy\, ds &\leq C \sum_{k=2}^\8 k^2 \int_0^\8 \int_{C_k} s^{-d/2} \exp\(-\frac{|z-x|^2}{4 s}\) \, dz \, ds\\
        &\leq C n^2 \int_{C_n} |z-x|^{2-d}\, dz + C k^2\sum_{2\leq k\neq n} \int_{C_k} |z-x|^{2-d} \, dz\\
        &= A_1 + A_2.
        }
    Observe that if $x\in D_n$ and $z\in C_n$, then $|x-z|\geq \tau/(2n)$. Therefore,
        \alx{
        A_1 &\leq C n^2 (\tau/n)^{2-d} n^{-d} = C \tau^{2-d} \leq c_0/2,
        }
    where the last inequality holds for $\tau$ big enough. Fix such $\tau$. In what follows we consider only $n\geq N_1$, such that $d(c_n, d_n)\leq 1/2$. For such $n$ and $k\neq n$ we have $|z-x| \geq c 2^{\max(n,k)}$ for $z\in C_k$ and $x\in D_n$. Thus,
        \alx{
        A_2 = \sum_{2\leq k<n}... + \sum_{k>n}... \leq & C \sum_{2\leq k<n} k^2 2^{n(2-d)} k^{-d} + C \sum_{k>n} k^2 2^{k(2-d)} k^{-d}\\
        &\leq Cn 2^{n(2-d)} + C 2^{n(2-d)}\leq c_0/2,
        }
    where the last estimate holds for $n\geq N_2$. The proof of \eqref{est4} finished by taking $N=\max(N_1,N_2)$.
}

Recall that $\qq^{[1]}$ consist of cubes of radii equal to 1 that satisfies $(G)$.
 We are now in position to prove that the spaces $H^1_{at}(\aa_{\qq^{[1]}})$ and $H^1_{at}(\aa_{\qq^{[1]},\w})$ are not equivalent as Banach spaces.

    \prop{main_ex}{
    There exist a sequence $a_n$ of $(\qq^{[1]},\w)-$atoms such that
        \eq{\label{atom_low_est}
        \norm{a_n}_{H^1_{at}(\aa_{\qq^{[1]}})} \geq c\ln n.
        }
    }

\begin{figure}
\includegraphics[width=250px]{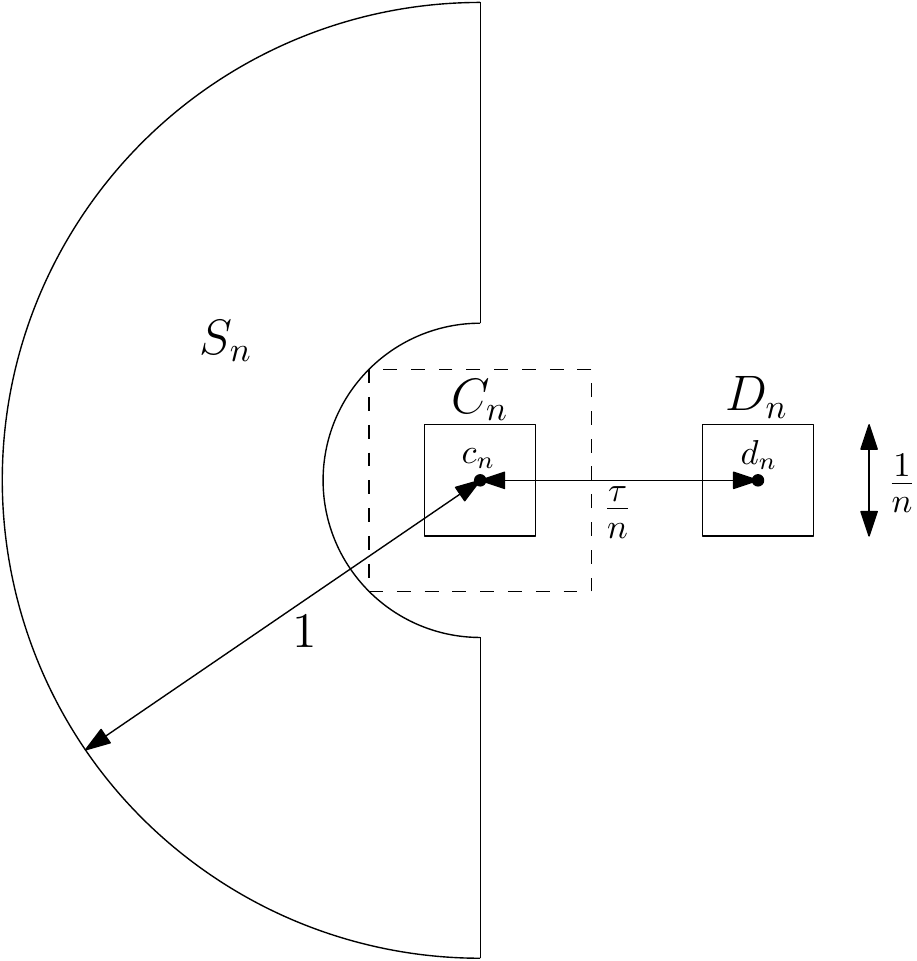}
\caption{The sets $C_n, D_n, S_n$.}
\label{figure1}
\end{figure}

    \pr{
    In this proof we use notation already introduced in Section \ref{sec6}. Let us denote $\w(S) = \int_S \w(x)\, dx$ and $\mu_n = \w(D_n) \w(C_n)^{-1}$. The atoms we are looking for are
        \eqx{
        a_n(x) = \zeta n^d \expr{\mu_n \chi_{C_n}(x) - \chi_{D_n}(x)},
        }
    where $\zeta >0$ is a constant that will be fixed in a moment.

    Let us check that $a_n$ are $(\qq^{[1]},\w)-$atoms. Obviously $\supp\, a_n \subseteq K_n := Q(c_n, (\tau + 1)/n)$. By the definition of $\mu_n$, $\int_{K_n} a_n(x) \w(x)\, dx = 0$. Recall that $|C_n| = |D_n|$, so by \eqref{omega_bounds} we get that $\mu_n \leq \delta^{-1}$. Moreover, by using Proposition \ref{p1}, for $n\geq N$,
        \eq{\label{mu_est}
        \mu_n \geq \frac{\inf\set{\w(x) \, : \, x\in D_n}}{\sup\set{\w(y) \, : \, y\in C_n}} = 1 + \frac{\inf\set{\w(x)-\w(y) \, : \, x\in D_n, \, y\in C_n}}{\sup\set{\w(y) \, : \, y\in C_n}} \geq 1+c_0.
        }
    What is left is to check the size condition. By choosing proper $\zeta>0$ we can write
        \eqx{
        \norm{a_n}_\8 \leq \zeta n^d \delta^{-1} \leq |K_n|^{-1},
        }
    so $a_n$ are indeed $(\qq,\w)-$atoms.

    Now we prove \eqref{atom_low_est}. For the collection $\qq^{[1]}$, the space $H^1_{at}(\aa_{\qq^{[1]}})$ is a classical local Hardy space. Equivalently, the norm can be given by a local maximal operator, see \eqref{Gold},
        \eqx{
        \norm{f}_{H^1_{at}(\qq^{[1]})} \simeq \norm{ \sup_{t\leq 1}\abs{\K{t}{0}f}}_{\LLL}.
        }
    Denote
        \eq{\label{def_Sn}
        S_n = \set{x\in \Rd \, : \, \sqrt{d}/n<\abs{x-c_n}<1, \ (x)_1<(c_n)_1 },
        }
    where $(x)_1$ is the first coordinate of $x\in \Rd$, see Figure \ref{figure1}. Obviously, $|S_n| \simeq C$. Assume now that $x\in S_n$ for some $n$. By \eqref{mu_est},
        \spx{
        \K{t}{0} a_n (x) = &\zeta n^d \int_\Rd (4\pi t)^{-d/2} \exp\expr{-\frac{|x-y|^2}{4t}}\expr{\mu_n \chi_{C_n}(y) - \chi_{D_n} (y)}\, dy\cr
        \geq& C n^d t^{-d/2} \int_\Rd \exp\expr{-\frac{|x-y|^2}{4t}}\expr{\chi_{C_n}(y) - \chi_{D_n} (y)}\, dy \cr
        &+ C n^d t^{-d/2} c_0 \int_{C_n} \exp\expr{-\frac{|x-y|^2}{4t}}\, dy\\
        =& A_1+A_2.
        }
    We claim that $A_1\geq 0$. Indeed, $D_n = C_n + (\tau/n)\mathbf{e_1}$ and for $x\in S_n$, $y_1\in C_n$ and $y_2= y_1+(\tau/n)\mathbf{e_1}$ we have $|y_1 -x| < |y_2 - x|$, c.f. \eqref{def_Sn}. We obtain that
        \eqx{
        A_1 = Cn^d t^{-d/2} \int_{C_n} \expr{ \exp\expr{-\frac{|x-y|^2}{4t}}-\exp \expr{ -\frac{|x-(y+\expr{\tau/n} \mathbf{e_1})|^2}{4t}}} \, dy \geq 0.
        }
        Now we deal with $A_2$. For $x\in S_n$ and $y\in C_n$ we have that $|x-y| \leq 2 |x-c_n|$. Thus,
        \eqx{
        A_2 \geq C t^{-d/2} \exp\expr{-\frac{|x-c_n|^2}{t}}.
        }
    Taking $t=|x-c_n|^2\leq 1$ we obtain that $\sup_{t\leq 1} A_2 \geq C |x-c_n|^{-d}$. The proof is finished by noticing that
        \eqx{
        \norm{\sup_{t\leq 1} \abs{\K{t}{0}a_n(x)}}_{L^1(S_n)} \geq C \int_{S_n} |x-c_n|^{-d}\, dx \geq C\,  \ln n,
        }
    where the last inequality is easily obtained by integrating in spherical coordinates.
    }

\section{Appendix}

    In the Appendix we consider a semigroup $\expr{\T{t}{}}_{t>0}$ that has positive integral kernel satisfying \eqref{up_gauss}. Obviously, all Schr\"odinger semigroups $\K{t}{U}$ with $0\leq U \in L^1_{loc}(\Rd)$ satisfy these assumptions.

    Our goal is to give a precise proof of the following natural estimate.

    \prop{prop_pointwise}{
    Assume that $f\in \LLL + L^\8(\Rd)$. For almost every $x\in \Rd$,
        \eqx{
        \lim_{t\to 0} \T{t}{f}(x) = f(x).
        }
    }

    \cor{cor_pointwise}{
    Let $0\leq U\in L^1_{loc}(\Rd)$ and $\tau>0$. Then
        \eqx{
        \norm{f}_\LLL \leq \norm{\Mm_\tau^U f}_\LLL.
        }
    }

    The proof of Proposition \ref{prop_pointwise} will be given at the end. We shall start with the following.

    \lem{lem_app}{
    Assume that $r>0$ is given. For a.e. $x\in \Rd$,
        \al{\label{kern_glob}
        \lim_{t\to 0} \int_{|x-y|>r} \t{t}{x,y}\, dy =0,\\
        \label{kern_loc}
        \lim_{t\to 0} \int_{|x-y|<r} \t{t}{x,y}\, dy =1.
        }
    }

    \pr{
    The equation \eqref{kern_glob} is a simple consequence of \eqref{up_gauss}. To prove \eqref{kern_loc} we shall use the fact that $\lim_{t\to 0} \T{t}f = f$, where the convergence is in $L^2(\Rd)$. From $L^2$ convergence we have a.e. convergence for a subsequence. Applying this to $f_n(x) = \chi_{Q(0, n)}(x)$, by a diagonal argument, we obtain a sequence $t_k>0$ that tends to zero, such that for a.e. $x\in\Rd$ we have
        \eq{\label{susseq}
        \lim_{k\to\8} \int_\Rd \t{t_k}{x,y}\, dy = 1.
        }
    Now, we are going to prove \eqref{susseq} for arbitrary sequence $s_j$ such that $\lim_{j\to \8} s_j =0$. Without loss of generality we can assume that $t_k$ is decreasing. For $j\in \Nn$, let $k_j$ be such that $t_{k_{j-1}}<s_j\leq t_{k_j}$ ($k_j=1$ when $s_j>t_{k_1}$). Then $t_{k_j} = s_j+ r_j$, where $\lim_{j\to \8} t_{k_j} = \lim_{j\to \8} r_j =0$. By \eqref{up_gauss} and the semigroup property,
        \spx{
        \int_\Rd \t{t_{k_j}}{x,y}\, dy = \int_\Rd \int_{\Rd} \t{s_j}{x,z} \t{r_j}{z,y}\, dz\, dy\leq \int_{\Rd} \t{s_j}{x,z} \leq 1.
        }
    Letting $j\to \8$, by \eqref{susseq}, we have that $\lim_{j\to \8} \int_\Rd \t{s_j}{x,z}\, dz =1$.
    }

    \pr{[Proof of Proposition \ref{prop_pointwise}]
    Assume that $f\in L^1(\Rd)+L^\8(\Rd) \subseteq L^1_{loc}(\Rd)$. By the Lebesgue differentiation theorem
        \eq{\label{ldt}
        \lim_{s\to 0} |Q(x,s)|^{-1} \int_{Q(x,s)} \abs{f(y)-f(x)} \, dy =0
        }
    for a.e. $x\in \Rd$. Assume that $x\in \Rd$ is such that \eqref{ldt}, \eqref{kern_glob} and \eqref{kern_loc} are satisfied for all rational $r>0$. The set of such points has full measure. For $\e>0$ fixed, we shall show that $\abs{\T{t} f(x) - f(x)}<C\e$ for $t$ small enough. Let $r>0$ be a fixed rational number such that for $s<r$ we have
        \eq{\label{ldt2}
        \int_{Q(x,s)}\abs{f(y)-f(x)}\, dy \leq \e |Q(x,s)|.
        }
    We can assume that $\sqrt{t}<r$. For such $t$, write
        \spx{
        \T{t}f(x) - f(x)=&f(x)\expr{\int_{|x-y|<r} \t{t}{x,y}\, dy-1} + \int_{|x-y|>r}\t{t}{x,y}f(y)\, dy\\
        & +\int_{|x-y|< \sqrt{t}}\t{t}{x,y}\expr{f(y)-f(x)}\, dy + \int_{\sqrt{t}<|x-y|< r} \t{t}{x,y}\expr{f(y)-f(x)}\, dy\\
         = &A_1+A_2+A_3+A_4.
        }
    By using \eqref{kern_loc}, we get that $A_1<\e$ for $t$ small enough. For the summand $A_2$ we consider two cases:
        \ite{
        \item if $f\in L^\8(\Rd)$, then $\abs{A_2}<\e$ for $t$ small enough by \eqref{kern_glob},
        \item if $f\in L^1(\Rd)$, then $\abs{A_2} \leq Ct^{-d/2} \exp\(-r^2/t\) \norm{f}_\LLL<\e$ for $t$ small enough.
        }
    By \eqref{up_gauss} and \eqref{ldt2}, for $t$ small enough,
        \eqx{
        A_3 \leq C t^{-d/2} \int_{|x-y|<\sqrt{t}} \abs{f(y)-f(x)}\, dy \leq C \e.
        }
    To estimate $A_4$ denote $N=\left\lceil \log_2 \frac{r}{\sqrt{t}}\right\rceil$, so that $r\leq \sqrt{t}2^N\leq 2 r$.
    Let $$R_n =\set{x\in \Rd \ : \  r2^{-n}<|x-y|<r2^{-n+1}}$$ for $n=1,...,N$. By \eqref{up_gauss} and \eqref{ldt2},
        \spx{
        A_4 &\leq C t^{-d/2} \sum_{n=1}^N \int_{R_n} \exp\expr{-\frac{|x-y|^2}{4t}}\abs{f(y)-f(x)}\, dy\\
        &\leq C t^{-d/2} \sum_{n=1}^N \exp\expr{-\frac{r2^{-n}}{ct}}\int_{R_n} \abs{f(y)-f(x)}\, dy\\
        &\leq C \e \sum_{n=1}^N \expr{\frac{r2^{-n}}{\sqrt{t}}}^{d} \exp\expr{-\frac{r2^{-n}}{c\sqrt{t}}}\leq C\e \frac{\sqrt{t}2^N}{r}\leq C\e.
        }
    }

{\bf Acknowledgments:} The author would like to thank Jacek Dziubański for discussions on the topic considered in the paper.

\bibliographystyle{amsplain}        
\bibliography{bib3}{}
\end{document}